\documentclass[11pt]{article}

\usepackage{graphicx}
\usepackage{amsmath,amsfonts,amssymb}
\usepackage{epstopdf}

\def\tto{\;{\lower 1pt \hbox{$\rightarrow$}}\kern -10pt
\hbox{\raise 2pt \hbox{$\rightarrow$}}\;}

\def\ve{\varepsilon}
\def\B{I\!\!B}
\def\h{\hfill\square}
\newcommand{\R}{\mathbb R}

\def\ox{\bar{x}}

\def\h{\hfill\square}

\def\O{\Omega}
\def\ph{\varphi}
\def\emp{\emptyset}

\newcounter{lk}

\begin{document}
\begin{center}
\vspace*{0.3in} {\bf A GENERALIZED SYLVESTER PROBLEM AND A
GENERALIZED FERMAT-TORRICELLI PROBLEM}\\[2ex]
Nguyen Mau Nam\footnote{Fariborz Maseeh Department of Mathematics and Statistics,
Portland State University, Portland, OR 97202, United States (email:
mau.nam.nguyen@pdx.edu). The research of Nguyen Mau Nam was partially
supported by the Simons Foundation under grant \#208785.}
 and Nguyen Hoang\footnote{Department of Mathematics, College of
Education, Hue University, Hue City, Vietnam (email:
nguyenhoanghue@gmail.com). The research of Nguyen Hoang was partially supported by the NAFOSTED, Vietnam, under grant \# 101.01-2011.26.}
\end{center}
\small{\bf Abstract:}  In this paper, we introduce and study the
following problem and its further generalizations: given two finite
collections of sets in a normed space, find a ball whose center lies
in a given constraint set with the smallest radius that encloses all
the sets in the first collection and intersects all the sets in the
second one. This problem can be considered as  a generalized version
of the \emph{Sylvester smallest enclosing circle problem} introduced
in the 19th century by Sylvester which asks for the circle of
smallest radius enclosing a given set of finite points in the plane.
We also consider a generalized version of the
\emph{Fermat-Torricelli problem}: given two finite collections of
sets in a normed space, find a point in a given constraint set that
minimizes the sum of the \emph{farthest distances} to the sets
in the first collection and \emph{ shortest distances (distances)} to the sets in the
second collection.

\medskip
\vspace*{0,05in} \noindent {\bf Key words.} The Sylvester smallest enclosing circle
problem, the Fermat-Torricelli problem, the smallest enclosing ball problem, the smallest intersecting ball problem

\newtheorem{Theorem}{Theorem}[section]
\newtheorem{Proposition}[Theorem]{Proposition}
\newtheorem{Remark}[Theorem]{Remark}
\newtheorem{Lemma}[Theorem]{Lemma}
\newtheorem{Corollary}[Theorem]{Corollary}
\newtheorem{Definition}[Theorem]{Definition}
\newtheorem{Example}[Theorem]{Example}
\renewcommand{\theequation}{\thesection.\arabic{equation}}
\normalsize

\section{Introduction and Problem Formulation}
\setcounter{equation}{0}

In the 19th century, the English mathematician Sylvester
(1814--1897) introduced the \emph{smallest enclosing circle
problem}: given a finite number of points in the plane, find the
smallest circle that encloses all of the points; see \cite{syl}. In
the 17th century, at the end of his book, \emph{Treatise on Maxima
and Minima}, the French mathematician Fermat (1601--1665) posed an
optimization problem which asks for a point that minimizes the sum
of the distances to three given points in the plane. This problem
was solved by the Italian mathematician and physicist Torricelli
(1608--1647) and is called the \emph{Fermat-Torricelli problem}.
These problems remain active as they are mathematically beautiful
and have meaningful real-world applications; see, e.g.,
\cite{AMS1,bms,chm,msw,frank,tan,wel} and the references
therein.\vspace*{0.05in}

Let $X$ be a normed space, and let $F\subset X$ be a closed,
bounded, convex set which contains the origin as an interior point.
For a point $x\in X$ and $r\geq 0$, the \emph{extended ball} with
center at $x$ and radius $r$ is defined by
\begin{equation*}
D_F(x; r):=x+rF.
\end{equation*}
It is obvious that when $F$ is the closed unit ball of $X$, the extended ball $D_F(x;r)$ reduces to the closed ball of radius $r$
centered at $x$.\vspace*{0.05in}

In the light of modern mathematics, we are going to study the
following problem called the \emph{generalized Sylvester problem}:
given two finite collections of sets and a constraint set in a
normed space, find an extended ball whose center lies in a given
constraint set with the smallest radius that encloses all the sets
in the first collection and intersects all the sets in the second
one. We also introduce and study the following generalized version
of the classical Fermat-Torricelli problem called the
\emph{generalized Fermat-Torricelli problem}: given two finite
collections of sets in a normed space, find a point in a given
constraint set that minimizes the sum of \emph{farthest distances}
the  to the sets in the first collection and \emph{shortest
distances (distances)} to the sets in the second
collection.\vspace*{0.05in}

Given a finite collection of nonempty, closed, bounded target sets
$\{\Omega_i: i\in I\}$ and another finite collection of nonempty,
closed target sets (not necessarily bounded) $\{\Theta_j: j\in J\}$,
and given a nonempty, closed constraint set $S$, the
\emph{generalized Sylvester problem} asks for a point $x\in S$ and
the smallest $r\geq 0$ such that
\begin{equation*}
\Omega_i\subset D_F(x;r)\mbox{ for all } i\in I, \mbox{ and
}D_F(x;r)\cap\Theta_j\neq\emptyset \mbox{ for all } j\in J.
\end{equation*}
In the case $I=\emptyset$, this problem reduces to the
\emph{smallest intersecting ball problem}, and in the case
$J=\emptyset$, it reduces to the \emph{smallest enclosing ball
problem}; see \cite{n1,n4}. It is also clear that when $X$ is the
Euclidean plane $\R^2$, all of the target sets under consideration
are singletons, and the constraint set $S=\R^2$, the generalized
Sylvester problem becomes the classical Sylvester enclosing circle
problem. Reducing to one of the simplest cases where three target
sets are three Euclidean balls in $\R^2$, and the constraint set
$S=\R^2$, we have shown in \cite{n4} that the solution of the
smallest intersecting ball problem has a close connection to the
solution of a particular case of the \emph{problem of Apollonius} on
tangent circles; see, e.g., \cite{gr}. It is interesting and not
hard to see that the generalized Sylvester considered in this paper
has a closed connection to the other cases of this celebrated
problem.

Following \cite{n1},  for a nonempty, closed, bounded set $Q \subset
X$, we define the \emph{maximal time function} to the target set $Q$
with the constant dynamic $F$ as follows:
\begin{equation}\label{maxima time function}
C_F(x; Q):=\inf \{ t\geq 0: Q\subset x+tF \}.
\end{equation}
The \emph{minimal time function} counterpart is  defined below as:
\begin{equation}\label{minimal time function}
T_F(x; Q):=\inf \{ t\geq 0: (x+tF)\cap Q\neq \emptyset \},
\end{equation}
where $Q$ needs not necessarily be bounded.\vspace*{0.05in}

Recall that the \emph{Minkowski function} associated with $F$ is defined by
\begin{equation*}
\rho_F(x):=\inf\{t\geq 0: x\in tF\}.
\end{equation*}
The maximal time function (\ref{maxima time function}) and the
minimal time function (\ref{minimal time function}) are Lipschitz
\\continuous and can be represented as
\begin{equation*}
C_F(x; Q)=\sup\{\rho_F(q-x): q\in Q\} \mbox{ and } T_F(x; Q)=\inf\{\rho_F(q-x): q\in Q\}.
\end{equation*}
Since $\rho_F(x)=\|x\|$ when $F$ is the closed unit ball of $X$, in this case, the maximal time function (\ref{maxima time function}) reduces to the \emph{farthest distance function}
\begin{equation*}
M(x;Q):=\sup\{\|x-q\|: q\in Q\},
\end{equation*}
and the minimal time function (\ref{minimal time function}) reduces to the \emph{distance function}
\begin{equation*}
d(x;Q):=\inf\{\|x-q\|: q\in Q\}.
\end{equation*}
In order to model the generalized Sylvester problem, we introduce the following function:
\begin{equation*}\label{c}
\mathcal{G}(x):=\max\{C_F(x;\Omega_i), T_F(x;\Theta_j): i\in I, j\in
J\},
\end{equation*}
and consider the following optimization problem
\begin{equation}\label{op}
\mbox{minimize }\mathcal{G}(x) \mbox{ subject to } x\in S.
\end{equation}
The \emph{generalized Fermat-Torricelli problem} can also be modeled as the following optimization problem:
\begin{equation}\label{hr}
\mbox{minimize } \mathcal{H}(x):=\sum_{i\in I}
C_F(x;\Omega_i)+\sum_{j\in J} T_F(x; \Theta_j)\mbox{ subject to }
x\in S.
\end{equation}
In this paper, we will mainly study the existence and uniqueness of
optimal solutions to the optimization problems (\ref{op}) and
(\ref{hr}) in general normed spaces. The existence results will be
studied in Section 2 and the uniqueness results will be studied in
Section 3. Our results generalize those obtained in \cite{n1,n3,n2,n4} and related references therein. \vspace*{0.05in}

We are going to use the following \emph{standing assumptions} throughout the paper:\\[1ex]
{\em $X$ is a normed space; $F$ is a closed, bounded, convex set
that contains $0$ as an interior point; $\Omega_i$ is nonempty,
closed, bounded for every $i\in I$; $\Theta_j$ is nonempty, closed
for every $j\in J$; S is a nonempty, closed set; and $I\cup J\neq
\emptyset$.}

\section{The Generalized Sylvester Problem and the Generalized Fermat-Torricelli Problem: the Existence of \\Optimal Solutions}
\setcounter{equation}{0}

In this section, we will study sufficient conditions that guarantee
the existence of optimal solutions to the optimization problems
(\ref{op}) and (\ref{hr}).

\begin{Lemma}\label{up} For $\alpha>0$, the level sets
\begin{equation*}
V_\alpha:=\{ x\in S: \mathcal{G}(x)<\alpha\}, \;W_\alpha:=\{ x\in S:
\mathcal{H}(x)<\alpha\}
\end{equation*}
and
\begin{equation*}
L_\alpha:=\{ x\in S: \mathcal{G}(x)\leq \alpha\}, \;K_\alpha:=\{
x\in S: \mathcal{H}(x)\leq \alpha\}
\end{equation*}
have the following estimates
\begin{equation*}
V_\alpha\subset S \cap \big[\cap_{i\in I}\cap_{\omega\in
\Omega_i}(\omega-\alpha F)\big ] \cap \big [ \cap_{j\in
J}(\Theta_j-\alpha F)\big]\subset L_\alpha,
\end{equation*}
and
\begin{equation*}
W_\alpha\subset S \cap \big[\cap_{i\in I}\cap_{\omega\in
\Omega_i}(\omega-\alpha F)\big ] \cap \big [ \cap_{j\in
J}(\Theta_j-\alpha F)\big]\subset K_{m\alpha},
\end{equation*}
where $m=|I|+|J|$.
\end{Lemma}
\noindent {\bf Proof: }Fix any $x\in V_\alpha$. It is obvious that $x\in S$ and the following hold for all $i\in I$ and for all $j\in J$:
\begin{equation*}
C_F(x;\Omega_i)<\alpha \mbox{ and } T_F(x;\Theta_j)<\alpha.
\end{equation*}
From the condition $C_F(x;\Omega_i)<\alpha$, one sees easily that
\begin{equation*}
\Omega_i\subset x+ \alpha F.
\end{equation*}
Thus, $$x\in \cap_{\omega\in \Omega_i}(\omega-\alpha F).$$
Similarly, from the condition $T_F(x;\Theta_j)<\alpha$, one finds
$0\leq t<\alpha$
\begin{equation*}
\Theta_j \cap (x+tF)\neq \emptyset.
\end{equation*}
Since $F$ is convex and $0\in F$, this implies
\begin{equation*}
x\in \Theta_j-tF\subset \Theta_j-\alpha F.
\end{equation*}
Now fix any $x\in S$ such that
\begin{equation*}
x\in \big[\cap_{i\in I}\cap_{\omega\in \Omega_i}(\omega-\alpha
F)\big ] \cap \big [ \cap_{j\in J}(\Theta_j-\alpha F)\big].
\end{equation*}
Then
\begin{equation*}
\Omega_i\subset x+\alpha F \mbox{ for all }i\in I, \mbox{ and }
(x+\alpha F)\cap \Theta_j\neq \emptyset \mbox{ for all }j\in J.
\end{equation*}
It follows that $C_F(x; \Omega_i)\leq \alpha$ for all $i\in I$, and
$T_F(x;\Theta_j)\leq \alpha$ for all $j\in J$. Thus,
$\mathcal{G}(x)\leq \alpha$, and hence $x\in L_\alpha$. The first
estimates have been proved. The second estimates can be proved in
the same way. The proof is now complete. $\h$

\begin{Lemma}\label{wlsc} Let $X$ be a reflexive Banach space.
Suppose $\Theta_j$ is weakly closed for every $j\in J$. Then
$\mathcal{G}$ and $\mathcal{H}$ are weakly sequentially lower
semicontinuous.
\end{Lemma}
{\bf Proof: }We will only show that $\mathcal{G}$ is weakly sequentially lower semicontinuous since the proof for $\mathcal{H}$ is similar. Fix any sequence $(x_k)$ that converges weakly to $\bar x$. We will show that
\begin{equation*}
\liminf_{k\to\infty}\mathcal{G}(x_k)\geq \mathcal{G}(\bar x).
\end{equation*}
Without loss of generality, suppose
\begin{equation*}
\liminf_{k\to\infty}\mathcal{G}(x_k)=\gamma\in \R.
\end{equation*}
Then there exists a subsequence of $(x_k)$ (without relabeling) such that
\begin{equation*}
\lim_{k\to\infty}\mathcal{G}(x_k)=\gamma.
\end{equation*}
For any $\ve>0$, there exists $k_0\in {\Bbb N}$ such that for any
$k\geq k_0$ and for any $i\in I$, $j\in J$, one has
\begin{equation*}
C_F(x_k;\Omega_i)<\gamma +\ve \mbox{ and }
T_F(x_k;\Theta_j)<\gamma+\ve.
\end{equation*}
It follows that for any $k\geq k_0$ and for any $i\in I$, $j\in
J$, the following hold:
\begin{equation*}
\Omega_i\subset x_k+(\gamma+\ve)F \mbox{ and }
(x_k+(\gamma+\ve)F)\cap \Theta_j\neq \emptyset.
\end{equation*}
The first inclusion implies
$$\Omega_i\subset \ox+(\gamma+\ve)F.$$ Since $X$ is reflexive, the set $F$ is a weakly sequentially
compact. Thus, the second condition implies $(\ox+(\gamma+\ve)F)\cap
\Theta_j\neq \emptyset$ under the assumption that $\Theta_j$ is
weakly closed for every $j\in J$. It follows that
\begin{equation*}
\mathcal{G}(\ox)=\max\{C_F(\ox;\Omega_i), T_F(\ox;\Theta_j): i\in I,
j\in J\}\leq \gamma+\ve.
\end{equation*}
Since $\ve>0$ is arbitrary, $\mathcal{G}(\ox)\leq \gamma$, and the
proof is now complete. $\h$ \vspace*{0.05in}

Let us define
\begin{equation*}
N(I,J,\alpha):=S \cap \big[\cap_{i\in I}\cap_{\omega\in
\Omega_i}(\omega-\alpha F)\big ] \cap \big [ \cap_{j\in
J}(\Theta_j-\alpha F)\big].
\end{equation*}
\begin{Proposition}\label{existence} The optimization problem {\rm(\ref{op})} has a nonempty optimal solution set under one of the following assumptions:\\[1ex]
{\rm (i)} There exists $\alpha>0$ such that $N(I,J, \alpha)$ is nonempty and precompact.\\
{\rm (ii)} $X$ is a reflexive Banach space; $S$ and $\Theta_j$ for $j\in J$ are weakly closed; and there exists $\alpha>0$ such that $N(I,J, \alpha)$ is nonempty and bounded.
\end{Proposition}
\noindent {\bf Proof: }We only need to prove the existence of an
optimal solution for (\ref{op}). The proof of the existence of an
optimal solution under (i) is straightforward since $\mathcal{G}$
is Lipschitz continuous; see, e.g., \cite{n1}. Let us give the
detail of the proof for the existence of an optimal solution under
(ii). Let $(x_k)\subset S$ be a minimizing sequence for the
optimization problem (\ref{op}). Since $N(I,J,\alpha)\neq
\emptyset$, and $N(I,J,\alpha)\subset L_\alpha$, by Lemma \ref{up},
\begin{equation*}
\inf\{\mathcal{G}(x): x\in S\}\leq \alpha.
\end{equation*}
In the case $\inf\{\mathcal{G}(x): x\in S\}= \alpha$, we see that any $x\in N(I,J,\alpha)$ is an optimal solution of the problem. In the other case,
\begin{equation*}
\mathcal{G}(x_k)<\alpha.
\end{equation*}
for all sufficiently large $k$. Thus, $x_k\in N(I,J,\alpha)$ for
such $k$, so $(x_k)$ is bounded. Since $X$ is reflexive and $S$ is
weakly closed, $(x_k)$ has a subsequence (without relabeling) that
converges weakly to $\bar x\in S$. Since the target sets $\Theta_j$
for $j\in J$ are also weakly closed, by Lemma \ref{wlsc},
\begin{equation*}
\mathcal{G}(\bar x)\leq \liminf_{k\to\infty}\mathcal{G}(x_k)=\inf\{\mathcal{G}(x): x\in S\}.
\end{equation*}
Therefore, $\bar x$ is an optimal solution of the problem. The proof is now complete. $\h$ \vspace*{0.05in}

The proof of the following proposition, with a slight difference in
formulation compared with Proposition \ref{existence}, is also
straightforward.

\begin{Proposition}\label{existence1} The optimization problem {\rm(\ref{hr})} has a nonempty optimal solution set under one of the following assumptions:\\[1ex]
{\rm (i)} There exists $\alpha>\inf\{\mathcal{H}(x): x\in S\}$ such that $N(I,J, \alpha)$ is precompact.\\
{\rm (ii)} $X$ is a reflexive Banach space; $S$ and $\Theta_j$ for $j\in J$ are weakly closed; and there exists $\alpha>\inf\{\mathcal{H}(x): x\in S\}$ such that $N(I,J, \alpha)$ is bounded.
\end{Proposition}

\begin{Theorem}\label{NC} The optimization problem {\rm(\ref{op})} and {\rm(\ref{hr})} have nonempty optimal solution sets under one of the following assumptions:\\[1ex]
{\rm (i)} The constraint set $S$ is compact.\\
{\rm (ii)} $I\neq\emptyset$ and $X$ is finite dimensional.\\
{\rm (iii)} $J\neq\emptyset$, at least one of the sets among $\{\Theta_j: j\in J\}$ is compact, and $X$ is finite dimensional.\\
{\rm (iv)} $X$ is a reflexive Banach space; $S$ and $\Theta_j$ for $j\in J$ are weakly closed; and at least one of them is bounded.
\end{Theorem}
{\bf Proof: }We will only prove the existence of an optimal solution
for (\ref{op}). Fix any \\$\alpha
> \inf\{\mathcal{G}(x): x\in S\}$. Then $N(I, J, \alpha)\neq
\emptyset$. We will show that assumption (i) of Proposition
\ref{existence} is satisfied under one of the assumptions: (i),
(ii), and (iii) in this theorem. Suppose (i) is satisfied. Then the
set $N(I, J, \alpha)$ is precompact since it is a subset of $S$. In
the case where (ii) is satisfied. Fix $i_0\in I$ and $\omega_0\in
\Omega_{i_0}$. Then
\begin{equation*}
N(I, J, \alpha)\subset \omega_0 -\alpha F.
\end{equation*}
This also implies that $N(I, J, \alpha)$ is precompact since
$\omega_0 -\alpha F$ is compact in this case. Now suppose that (iii)
is satisfied. Choose $j_0\in J$ such that $\Theta_{j_0}$ is compact.
Then, again,  $N(I, J, \alpha)$ is precompact since $N(I, J,
\alpha)\subset \Theta_{j_0}-\alpha F$ and the latter is compact in
this case. To finish the proof, we will show that condition (ii) of
Proposition \ref{existence} is satisfied under condition (iv) of
this theorem. If $J=\emptyset$, then $S$ is bounded, so $N(I, J,
\alpha)\subset S$ is bounded. In the case $J\neq \emptyset$, it is
also easy to see that $N(I, J, \alpha)$ is bounded. The proof is now
complete. $\h$ \vspace*{0.05in}

Finally, we will study the relationship between the optimization problem (\ref{op}) and the generalized Sylvester problem in the proposition below.

\begin{Proposition} Let $X$ be a reflexive Banach space. Suppose that $\Theta_j$ is weakly closed for every $j\in J$.
Then $\ox$ is an optimal solution of the generalized Sylvester problem with radius $r$ if and only if $\ox$ is an optimal
solution of the optimization problem {\rm(\ref{op})} with
$r=\mathcal{G}(\ox).$
\end{Proposition}
{\bf Proof:} Let $\ox\in S$ be an optimal solution of the generalized Sylvester problem with the smallest radius $r\geq 0$. Then
\begin{equation*}
\Omega_i \subset \ox +rF \mbox{ for all }\ i \in I, \mbox{ and
}\Theta_j \cap (\ox +rF)\ne\emptyset \mbox{ for all } j\in J.
\end{equation*}
Moreover, for any $x\in S$ and $t\ge 0$ such that
\begin{equation*}
\Omega_i \subset x +tF \mbox{ for all } i\in I, \mbox{ and }\Theta_j
\cap (x +tF)\ne \emp \mbox{ for all } j\in J,
\end{equation*}
one has $r\le t$. We will show that $\mathcal{G}(\ox) =r$ and $\mathcal G(\ox) \le \mathcal G(x) \mbox{ for all }x\in S.$

Since
\begin{equation*}
\Omega_i \subset \ox +rF \mbox{ for all }i\in I, \mbox{ and }
\Theta_j\cap (\ox +rF)\ne \emp \mbox{ for all }j\in J,
\end{equation*}
one has $\mathcal{G}(\ox)\le r$.

Since $0\in F$ and $F$ is convex, $t_1F\subset  t_2F$ whenever $0\le t_1\le t_2$. If $\mathcal{G}(\ox)<r$, take $r'>0$ such that
$\mathcal{G}(\ox) <r'<r$. Then
$$\Omega_i\subset \ox +r'F \mbox{ for all }i\in I, \mbox{ and } (\ox +r'F)\cap \Theta_j \ne \emp \mbox{ for all }j\in J.$$
This contradicts to the definition of an optimal solution of the generalized Sylvester problem.

Now fix any $x\in S$ and define $r':=\mathcal{G}(x)$. It is not hard to see that
\begin{equation*}
\Omega_i\subset \ox +r'F \mbox{ for all }i\in I, \mbox{ and } (\ox
+r'F)\cap \Theta_j \ne \emp \mbox{ for all }j\in J.
\end{equation*}
Thus, $\mathcal{G}(\ox)=r\leq r'=\mathcal{G}(x)$. \vspace*{0.05in}

Conversely, suppose that $\ox \in S$ is a solution of the optimal problem (\ref{op}). Then
\begin{equation*}
\mathcal{G}(\ox)=\min_{x\in S}\mathcal{G}(x)=r.
\end{equation*}
Thus,
$$C_F(\ox, \Omega_i)\le r \mbox{ for all }i\in I, \mbox{ and }T_F(\ox,\Theta_j)\le r \mbox{ for all }j\in J.$$
Consequently, $$\Omega_i\subset \ox +rF \mbox{ for all }i\in
I,\mbox{ and }\Theta_j\cap (\ox+rF)\neq\emptyset\mbox{ for all }j\in J.$$ Now,
take $x\in S$ and $t\ge 0$ such that
$$\Omega_i \subset x +tF \mbox{ for all }i\in I \mbox{ and }\Theta_j\cap (x +tF)\ne \emp \mbox{ for all }j\in J.$$
Let $r' =\mathcal{G}(x)$. Then $r\le r'\le t.$ This means  $\ox$ is an optimal solution of the generalized Sylvester problem. $\h$

\section{The Generalized Sylvester Problem and the Generalized Fermat-Torricelli Problem: the Uniqueness of \\Optimal Solutions}
\setcounter{equation}{0}

We are first going to study the uniqueness of an optimal solution to the optimization problem (\ref{op}). Recall that a set $C$ is called \emph{convex} if for every $x, y\in C$, one has
\begin{equation*}
[x,y]:=\{tx+(1-t)y: t\in [0,1]\}\subset C.
\end{equation*}
The set $C$ is called \emph{strictly convex} if for every $x,y\in
C$, $x\neq y$, and for every $t\in (0,1)$, one has
\begin{equation*}
tx+(1-t)y\in \mbox{int } C.
\end{equation*}
A function $f$ is called \emph{convex} on a convex set
$S$ if for every $x,y\in S$ and for every $t\in (0,1)$, one has
\begin{equation*}
f(tx+(1-t)y)\leq tf(x)+(1-t)f(y).
\end{equation*}
If this inequality becomes strict for all $x, y\in S$, $x\neq y$,
and for every $t\in (0,1)$, the function is called \emph{strictly
convex}.\vspace*{0.05in}

The following lemma is useful in the sequel. It is also of independent interest.

\begin{Lemma}\label{strict convexity} Let $S$ be convex. Suppose  $g$ is a nonnegative, convex function on  $S$. Then the function defined by
 $$h(x):=(g(x))^2$$
 is strictly convex if and only if $g$ is not constant on any line segment $[a,b]\subset S$, where $a\neq b$.
\end{Lemma}
\noindent {\bf Proof:} Suppose $h$ is strictly convex on $S$. On the contrary, suppose that $g$ is constant on a line segment $[a, b]$, where $a\neq b$. Then $h$ is also constant on this line segment, which is a contradiction.

Conversely, suppose $g$ is not constant on any line segment $[a,b]\subset S$ , where $a\neq b$. For any $t\in (0,1)$ and for any $x,y\in S$, one has
\begin{align*}
h(tx+(1-t)y)&=(g(tx+(1-t)y))^2\le (tg(x) +(1-t)g(y))^2\\
&= t^2(g(x))^2 +(1-t)^2(g(y))^2 +2t(1-t)g(x)g(y)\\
&\le t^2(g(x))^2 +(1-t)^2(g(y))^2 +t(1-t) ((g(x))^2 +(g(y))^2)\\
&= t(g(x))^2 +(1-t) (g(y))^2= th(x) +(1-t)h(y).
\end{align*}
Thus, $h$ is a convex function on $S$. We will show that it is
strictly convex on $S$. Suppose by contradiction that there exist
$t\in (0,1)$ and $x,y\in S$ such that
\begin{equation*}
h(z)=th(x)+(1-t)h(y),
\end{equation*}
where $z:=tx+(1-t)y.$ Then
\begin{equation*}
(g(z))^2=t(g(x))^2+(1-t)(g(y))^2.
\end{equation*}
Since $g(z)\le tg(x)+(1-t)g(y)$, one has
\begin{equation*}
(g(z))^2\le t^2(g(x))^2 +(1-t)^2(g(y))^2 +2t(1-t)g(x)g(y).
\end{equation*}
This implies
\begin{equation*}
t(g(x))^2+(1-t)(g(y))^2 \le  t^2(g(x))^2 +(1-t)^2(g(y))^2 +2t(1-t) g(x)g(y).
\end{equation*}
Thus, $(g(x) -g(y))^2\le 0$, and hence $g(x)=g(y).$ We have proved that $h(x) =h(z) =h(y)$, where $z\in (x,y)$. We will get a contradiction by showing that $h$ is constant on the line segment $[z,y].$ Indeed, fix any $u\in (z,y)$. Then
$$h(u)\le \nu h(z) +(1-\nu)h(y) =h(z)\mbox{ for some }\nu\in (0,1).$$
On the other hand, since $z$ lies in between $x$ and $u$, one has
\begin{equation*}
h(z)\le \mu h(x) +(1-\mu)h(u)\le \mu h(z) +(1-\mu)h(z)=h(z) \mbox{ for some } \mu\in (0,1).
\end{equation*}
Thus, $h(z)=\mu h(x) +(1-\mu)h(u)=\mu h(z)+(1-\mu)h(u)$, and hence $h(u)=h(z)$. This contradicts to the assumption that $g(u)=\sqrt{h(u)}$ is not constant on any line segment $[a,b]$, where $a\neq b$. The proof is now complete. $\h$ \vspace*{0.05in}

Given a nonempty, closed, bounded subset $\Omega$ of $X$ and a point
$x\in X$, the \emph{farthest projection} from $x$ to $\Omega$ is
defined by
\begin{equation*}
\mathcal{P}_F(x;\Omega):=\{\omega\in \Omega:
\rho_F(\omega-x)=C_F(x;\Omega)\}.
\end{equation*}

\begin{Proposition}\label{constant C}
Suppose $F$ is strictly convex. Let $\Omega$ be a nonempty, closed,
bounded subset of $X$ such that $\mathcal{P}_F(x;\Omega)\neq
\emptyset$ for all $x\in S$, where $S$ is convex. Then the function
$h(x):=C_F(x;\Omega)$ is convex and not constant on any straight
line segment $[a,b]\subset S$, where $a\neq b$.
\end{Proposition}
\noindent {\bf Proof:} The function $h$ is obviously convex on $S$
since it is the supremum of a family of convex functions by the
representation
\begin{equation*}
C_F(x;\Omega)=\sup\{\rho_F(\omega-x): \omega\in \Omega\}.
\end{equation*}
Suppose by contradiction that $h(x)=r$ for all $x\in [a,b]\subset S$
for some line segment $[a,b]$, where $a\neq b$. Let $\omega\in
\mathcal{P}_F(\dfrac{a+b}{2};\Omega)$. Then
$$r = h(\frac{a+b}{2}) =\rho_F(\omega-\dfrac{a+b}{2}).$$
It follows that
\begin{equation*}
r=\rho_F(\dfrac{\omega-a}{2}+\dfrac{\omega-b}{2})\leq
\dfrac{1}{2}\rho_F(\omega-a) +\dfrac{1}{2}\rho_F(\omega-b)\leq
\dfrac{1}{2}(C_F(a;\Omega)+C_F(b;\Omega))=r.
\end{equation*}
Since $r=C_F(a;\Omega)\geq \rho_F(\omega-a)$ and
$r=C_F(b;\Omega)\geq \rho_F(\omega-b)$,
\begin{equation*}
\rho_F(\omega-a)=\rho_F(\omega-b)=r.
\end{equation*}
Since $r>0$, this implies
\begin{equation*}
\rho_F(\dfrac{\omega-a}{r})=\rho_F(\dfrac{\omega-b}{r})=1.
\end{equation*}
Thus, $\dfrac{\omega-a}{r}\in F$ and $\dfrac{\omega-b}{r}\in F$. Since $F$ is strictly convex,
\begin{equation*}
\dfrac{1}{2}(\dfrac{\omega-a}{r}) +\dfrac{1}{2}(\dfrac{\omega-b}{r})=\dfrac{1}{r}(\omega-\dfrac{a+b}{2}) \in \mbox{int }F.
\end{equation*}
This implies
\begin{equation*}
\rho_F(\dfrac{1}{r}(\omega-\dfrac{a+b}{2}))<1, \mbox{ and hence }\rho_F(\omega-\dfrac{a+b}{2})<r,
\end{equation*}
which is a contradiction. The proof is now complete. $\h$ \vspace*{0.05in}

For a nonempty, closed subset $\Theta$ of $X$ and $x\in X$, the
\emph{projection} from $x$ to $\Theta$ is defined by
\begin{equation*}
\Pi_F(x; \Theta):=\{ u\in \Theta: \rho_F(u-x)=T_F(x;\Theta)\}.
\end{equation*}
It is not hard to see that $\Pi_F(x; \Theta)\neq \emptyset$ for all $x\in X$ under one of the following conditions:\\
(i) $\Theta$ is compact.\\
(ii) $X$ is finite dimensional and $ \Theta$ is closed.\\
(iii) $X$ is reflexive and $\Theta$ is weakly closed.

\begin{Proposition}\label{constant T} Suppose $F$ is strictly convex.
Let $\Theta$ be a nonempty, closed, strictly convex set of $X.$
Suppose further that the set $\Pi_F(x,\Theta)$ is nonempty for all
$x\in S$, where $S$ is a convex set. Then the function
$g(x)=T_F(x,\Theta)$ is convex and not constant on any line segment
$[a,b]\subset S$ such that $a\neq b$ and $[a,b]\cap \Theta=\emp.$
\end{Proposition}
\noindent {\bf Proof:} It is easy to see that $g(x)$ is convex on
$S$. Observe that if $u\in \Pi_F(x;\Theta)$, where $x\notin\Theta$,
then $u\in\mbox{\rm bd } \Theta$. Indeed, if $u\in
\text{int}\;\Theta$, then there exists $\ve>0$ such that
\begin{equation*}
\B(u,\ve)\subset \Theta.
\end{equation*}
Define
\begin{equation*}
z:=u+\ve\dfrac{x-u}{\|x-u\|}\in \B(u,\ve).
\end{equation*}
Then
\begin{equation*}
\rho_F(z-x)=\rho_F(u+\ve\dfrac{x-u}{\|x-u\|}-x)\| =(1-\frac{\ve}{\|x-u\|})\rho_F(u-x) <\rho_F(u-x)=T_F(x;\Theta),
\end{equation*}
for $\ve$ sufficiently small, which is a contradiction.

On the contrary, suppose that there exists $[a,b]\subset S$, $a\neq
b$, such that $[a,b]\cap\Theta\neq\emptyset$, and $T_F(x;\Theta)=r$
for all $x\in [a,b].$ Choose $u\in \Pi_F(a;\Theta)$ and $v\in
\Pi_F(b; \Theta)$. Then
\begin{equation*}
\rho_F(u-a)=\rho_F(v-b)=r.
\end{equation*}
We will first show that $u\neq v$. Indeed, if $u=v$, then
\begin{equation*}
\rho_F(u-a)=\rho_F(u-b)=r.
\end{equation*}
Following the proof of Lemma \ref{constant C}, one has
\begin{equation*}
\rho_F(u-\dfrac{a+b}{2})<r,
\end{equation*}
and hence
$$T_F(\dfrac{a+b}{2}; \Theta)\leq \rho_F(u-\dfrac{a+b}{2})<r,$$
which is not the case.

For a fixed $t\in (0,1)$, one has
\begin{align*}
r=T_F(ta+(1-t)b;\Theta)&\leq \rho_F(tu+(1-t)v-(ta+(1-t)b))\\
&\leq t\rho_F(u-a) +(1-t)\rho_F(v-b)=r.
\end{align*}
This implies $tu+(1-t)v\in \Pi_F(ta+(1-t)b; \Theta)$. Thus, $tu+(1-t)v\in \mbox{bd }\Theta$. This contradicts the strict convexity of $\Theta$. The proof is now complete. $\h$

\begin{Lemma}\label{constant} Suppose that $h_i$ for $i=1,\ldots,m$, $m\geq 1$, are nonnegative, continuous, convex functions on $S$, where $S$ is convex. Define
\begin{equation*}
\phi(x): =\max\{h_1(x),\dots, h_m(x)\}.
\end{equation*}
 Suppose that $\phi(x) =r>0$ for all $x\in [a,b]$ for some line segment $[a,b]\subset S$, where $a\neq b$. Then there exists a line segment $[\alpha,\beta]\subset [a,b]$, $\alpha\neq\beta$,  and $i_0\in\{1,\ldots, m\}$ such that $$h_{i_0}(x)=r \mbox{ for all }x\in [\alpha, \beta].$$
\end{Lemma}
\noindent {\bf Proof:} The conclusion is obvious for $m=1$. Suppose
that $$\phi (x)=\max\{h_1(x), h_2(x)\}.$$ The conclusion is
obviously true if $h_1(x)=r$ for all $x\in [a,b].$ Otherwise, there
exists $x_0\in [a,b]$ such that $h_1(x_0)< r.$ Then there exists a
subinterval $[\alpha, \beta ]\subset [a,b]$, $\alpha\neq\beta$, such
that
\begin{equation*}
h_1(x)<r \mbox{ for all } x\in [\alpha, \beta].
\end{equation*}
Therefore, $h_2(x) =r$ on this subinterval. Suppose that the conclusion holds for a positive integer $m$. Let $\phi (x)=\max\{h_1(x), \dots, h_m(x), h_{m+1}(x)\}$. Then $\phi(x)=\max \{h_1(x), k_1(x)\}$ where $k_1(x):=\max\{h_2(x),\dots, h_{m+1}(x)\}.$ The conclusion follows from the case $m=2$ and the induction assumption. The proof is now complete. $\h$ \vspace*{0.05in} \vspace*{0.05in}

Now we are ready to prove our main theorem on a necessary and sufficient conditions
for the optimization problem (\ref{op}) to have at most one optimal
solution. To obtain sufficient conditions for the uniqueness of an
optimal solution to this problem, we only need to combine this theorem with the
results from Theorem \ref{NC}.

\begin{Theorem}\label{uniqueness} Let $F$ and $\Theta_j$ for $j\in J$ be strictly convex, and let $S$ be convex.
Suppose that for every $x\in S,$  the projection sets
$\mathcal{P}_F(x;\Omega_i)$ and $\Pi_F(x,\Theta_j)$ are not empty for
all $i\in I$ and $j\in J$. Then the optimization problem {\rm
(\ref{op})} has at most one optimal solution if and only if one of the following conditions is satisfied:\\[1ex]
{\rm (1)} The index set $I$ is empty and $\cap_{j\in J}[\Theta_j\cap S]$ contains at most one point.\\
{\rm (2)} The index set $I$ is nonempty.
\end{Theorem}
{\bf Proof:} Define
\begin{equation*}
\mathcal{C}_1(x):=\max\{C_F(x;\Omega_i): i\in I\}
\end{equation*}
and
\begin{equation*}
\mathcal{T}_1(x):=\max\{T_F(x; \Theta_j): j\in J\}.
\end{equation*}
Then
\begin{equation*}
\mathcal{G}(x)=\max\{\mathcal{C}_1(x), \mathcal{T}_1(x)\}.
\end{equation*}
We also defined
\begin{equation*}
\mathcal{K}(x):=(\mathcal{G}(x))^2=\max\{(\mathcal{C}_1(x))^2, (\mathcal{T}_1(x))^2\}.
\end{equation*}
Consider the optimization problem
\begin{equation}\label{sq}
\mbox{minimize } \mathcal{K}(x) \mbox{ subject to }x\in S.
\end{equation}
It is obvious that $\ox\in S$ is an optimal solution of problem (\ref{op}) if and only if it is an optimal solution to problem (\ref{sq}). We are going to prove that $\mathcal{K}$ is strictly convex on $S$ under (1) or (2).

Suppose first that (1) is satisfied. In this case, it suffices to show that $\ph(x):=(\mathcal{T}_1(x))^2$ is strictly convex on the set $S$. By contradiction, suppose $\ph$ is not strictly convex on $S$. By
Lemma \ref{strict convexity}, there exists a line segment
$[a,b]\subset S$, $a\neq b$, and $r\geq 0$ such that
\begin{equation*}
\mathcal{T}_1(x)=\max\{T_F(x;\Theta_j): j\in J\}=r \mbox{ for all
}x\in [a, b].
\end{equation*}
It is clear that $r>0$, since otherwise, $[a,b]\subset \cap_{j\in
J}\Theta_j$, which is a contradiction since $\cap_{j\in J}[\Theta_j\cap S]$ contains at most one point. By Lemma
\ref{constant}, there exists a line segment $[\alpha, \beta]\subset
[a,b]$, $\alpha\neq\beta$, and $j_0\in J$ such that
\begin{equation*}
T_F(x;\Theta_{j_0})=r \mbox{ for all }x\in [\alpha, \beta].
\end{equation*}
Since $r>0$, $x\notin \Theta_{j_0}$ for all $x\in [\alpha, \beta]$. Thus, $[\alpha, \beta]\cap \Theta_{j_0}=\emptyset$. This is a contradiction to Proposition \ref{constant T}.

Now, let us assume that (2) is satisfied. We will show that $\mathcal{K}$ is also strictly convex on $S$ in this case. Again, by contradiction, suppose $\mathcal{K}$ is not strictly convex on $S$. By
Lemma \ref{strict convexity}, there exists a line segment
$[a,b]\subset S$, $a\neq b$, and $r\geq 0$ such that
\begin{equation*}
\mathcal{G}(x)=\max\{\mathcal{C}_1(x), \mathcal{T}_1(x)\}=r \mbox{ for all
}x\in [a, b].
\end{equation*}
In the case where $r=0$, one has $\mathcal{C}_1(x)=0$ for all $x\in [a, b]$. Thus, for any $i\in I$, one has $C_F(x;\O_i)=0$ for all $x\in [a, b]$. In this case $\O_i$ must be a singleton and that contains $x$ for all $x\in [a, b]$, which is not the case. In the case where $r>0$, by Lemma
\ref{constant}, one of the following holds:\\[1ex]
(a) There exist a line segment $[\alpha, \beta]\subset
[a,b]$, $\alpha\neq\beta$, and $j_0\in J$ such that
\begin{equation*}
T_F(x;\Theta_{j_0})=r \mbox{ for all }x\in [\alpha, \beta].
\end{equation*}
(b) There exist a line segment $[\alpha, \beta]\subset
[a,b]$, $\alpha\neq\beta$, and $i_0\in I$ such that
\begin{equation*}
C_F(x;\Omega_{i_0})=r \mbox{ for all }x\in [\alpha, \beta].
\end{equation*}
If (a) holds, then we arrive at a contradiction to Proposition \ref{constant T} in the same way as the previous proof. In the case (b) holds, we also arrive at a contradiction to Proposition \ref{constant C}. We have shown that the function $\mathcal{K}$ is strictly convex on $S$, and hence problem (\ref{op}) has at most one optimal solution.

Let us now show that if problem (\ref{op}) has at most one solution, then either (1) or (2) is satisfied. Suppose by contradiction that both conditions are not satisfied (notice that (1) and (2) cannot occur simultaneously). Then index set $I=\emptyset$ and $\cap_{j\in J}[\Theta_j\cap S]$ contains more than one points. It clear that $\cap_{j\in J}[\Theta_j\cap S]$ is the solution set of the problem in this case with the optimal value zero. Thus, the problem has more than one solution, which is a contradiction. The proof is now complete. $\h$ \vspace*{0.05in}

Theorem \ref{uniqueness} generalizes the sufficient conditions
for the uniqueness of optimal solutions for the smallest enclosing
ball problem and the smallest intersecting ball problem given in
\cite{n1}.\vspace*{0.05in}

The next two corollaries follow directly from Theorem
\ref{uniqueness}.

\begin{Corollary} Let $F$ be a strictly convex set, let $S$ be convex, and let $I\neq\emptyset$, $J\neq\emptyset$. \\[1ex]
{\rm (i)} Suppose that for every $x\in S,$  the projection set
$\mathcal{P}_F(x;\Omega_i)$  is nonempty for all $i\in I$. Then the
optimization problem
\begin{equation*}
\mbox{\rm minimize }\mathcal{C}_1(x)=\max\{C_F(x;\Omega_i): i\in I\}
\mbox{\rm \;subject to }x\in S
\end{equation*}
has at most one optimal solution.\\[1ex]
{\rm (ii)} Suppose that $\Theta_j$ is
strictly convex for every $j\in J$, and for every $x\in S,$  the projection set
$\Pi_F(x,\Theta_j)$ is nonempty for all $j\in J$. Then the
optimization problem
\begin{equation*}
\mbox{\rm minimize }\mathcal{T}_1(x)=\max\{T_F(x;\Theta_j): j\in
J\}
\mbox{\rm\; subject to }x\in S
\end{equation*}
has at most one optimal solution if and only if $\cap_{j\in J}[\Theta_j\cap S]$ contains at most one point.
\end{Corollary}
It seems that our results are new even when reducing to the simple case below.
\begin{Corollary}  Suppose $F$ is strictly convex and $S$ is convex. For any finite collection of points $\{a_i: i=1,\ldots,m\}$, $m\geq 1$, the optimization
\begin{equation*}
\mbox{\rm minimize }\max\{\rho_F(a_i-x): i=1,\ldots,m\} \mbox{\rm \; subject to }x\in S
\end{equation*}
has at most one optimal solution. This problem has a unique optimal solution if we assume additionally that $X$ is reflexive or $S$ is compact.
\end{Corollary}

Recall that $X$ is called strictly convex if  the closed unit ball
of $X$ is a strictly convex set. It is well-known that $X$ is
strictly convex if and only if the following implication holds:
\begin{equation*}
[x\ne  y \mbox{ and }\|x\|=\|y\|=1]\Rightarrow \| x + y \| < 2.
\end{equation*}
There are several examples of strictly convex normed spaces such as
Hilbert spaces and $L^p$ spaces for $p>1$.\vspace*{0.05in}

The following corollary follows directly from Theorem \ref{uniqueness}.

\begin{Corollary} Let $\Theta_j$ be strictly convex for every $j\in J$ and let $S$ be convex.
Suppose that $X$ is strictly convex and $F$ is the closed unit ball
of $X$. Suppose further that for every $x\in S,$  the projection
sets $\mathcal{P}_F(x;\Omega_i)$ and $\Pi(x,\Theta_j)$ are nonempty
for all $i\in I$ and $j\in J$. Then the optimization problem {\rm
(\ref{op})} has at most one optimal solution if and only if $I\neq\emptyset$ or $\cap_{j\in J}[\Theta_j\cap S]$ contains at most one point.
\end{Corollary}

We are now going to study the uniqueness of an optimal solution to
the optimization problem (\ref{hr}). The following lemma will be
important for the study. It generalizes a familiar property of the
norm on a strictly convex normed space.

\begin{Lemma}\label{ln} Suppose $F$ is strictly convex. If $x\neq 0$ and $y\neq 0$, then
\begin{equation}\label{ta}
\rho_F(x+y)=\rho_F(x)+\rho_F(y)
\end{equation}
if and only if $x=\lambda y$ for some $\lambda >0$.
\end{Lemma}
\noindent {\bf Proof: }Let $\alpha:=\rho_F(x)$ and $\beta:=\rho_F(y)$. Since $F$ is bounded, $\alpha>0$ and $\beta>0$. From (\ref{ta}), one has
\begin{equation*}
\rho_F(\dfrac{x+y}{\alpha+\beta})=1.
\end{equation*}
This implies
\begin{equation*}
\rho_F(\dfrac{x}{\alpha}\dfrac{\alpha}{\alpha+\beta}+\dfrac{y}{\beta}\dfrac{\beta}{\alpha+\beta})=1.
\end{equation*}
Thus,
\begin{equation*}
\dfrac{x}{\alpha}\dfrac{\alpha}{\alpha+\beta}+\dfrac{y}{\beta}\dfrac{\beta}{\alpha+\beta}\in \mbox{bd }F.
\end{equation*}
Since $\dfrac{x}{\alpha}\in F$, $\dfrac{y}{\beta}\in F$, and $\dfrac{\alpha}{\alpha+\beta}\in (0,1)$, the strict convexity of $F$ implies $$\dfrac{x}{\alpha}=\dfrac{y}{\beta}.$$
It follows that $x=\lambda y$ for $\lambda:=\dfrac{\rho_F(x)}{\rho_F(y)}$. The opposite implication is obvious. The proof is now complete. $\h$

\begin{Lemma} Suppose $F$ is strictly convex and $\Theta$ is convex. Then for any $x\in X$, the set $\Pi_F(x;\Theta)$
cannot contain more than one points.
\end{Lemma}
\noindent {\bf Proof:} We only need to consider the case where $x\notin \Theta$. Suppose by contradiction that there exist $u_1, u_2\in \Pi_F(x;\Theta)$ and $u_1\neq u_2$. Then
\begin{equation*}
T_F(x;\Theta)=\rho_F(u_1-x)=\rho_F(u_2-x)=r>0.
\end{equation*}
This implies $\dfrac{u_1-x}{r}\in F$ and $\dfrac{u_2-x}{r}\in F$. Since $F$ is strictly convex,
\begin{equation*}
\dfrac{1}{2}\dfrac{u_1-x}{r} +\dfrac{1}{2}\dfrac{u_2-x}{r}=\dfrac{1}{r}(\dfrac{u_1+u_2}{2}-x)\in \mbox{int }F.
\end{equation*}
It follows that $\rho_F(u-x)<r=T_F(x; \Theta)$, where
$u=\dfrac{u_1+u_2}{2}\in \Theta$. This is a contradiction.
$\h$\vspace*{0.05in}

In what follows, we identify the projection $\Pi_F(x;\Theta)$ with
its unique element when $F$ and $\Theta$ are strictly
convex.\vspace*{0.05in}

For two different points $x$ and $y$ in $X$, define
\begin{equation*}
L(x,y):=\{tx+(1-t)y: t\in {\Bbb R}\}.
\end{equation*}

\begin{Proposition} Let $F$ be strictly convex, let $I\neq\emptyset$, and let $S$ be convex.
Suppose that for any $x, y\in S$, $x\neq y$, there exists $i\in I$
such that
\begin{equation*}
L(x,y)\cap \Omega_i=\emptyset.
\end{equation*}
Under the assumption that for every $x\in X$, the projection set
$\mathcal{P}_F(x; \Omega_i)\neq \emptyset$ for all $i\in I$, the
function
\begin{equation*}
\mathcal{C}_2(x):=\sum_{i\in I} C_F(x;\Omega_i)
\end{equation*}
is strictly convex on $S$.
\end{Proposition}
\noindent {\bf Proof: }Suppose by contradiction that $\mathcal{C}_2$
is not strictly convex on $S$. Then there exist $x, y\in S$, $x\neq
y$, and $t\in (0,1)$ such that
\begin{equation*}
\mathcal{C}_2(tx+(1-t)y)=t\mathcal{C}_2(x)+(1-t)\mathcal{C}_2(y).
\end{equation*}
This implies
\begin{equation}\label{quali}
C_F(tx+(1-t)y;\Omega_i)=tC_F(x;\Omega_i)+(1-t)C_F(y;\Omega_i) \mbox{
for all }i \in I.
\end{equation}
Choose $i_0\in I$ such that
\begin{equation*}
L(x,y)\cap \Omega_{i_0}=\emptyset.
\end{equation*}
For any $\omega\in \mathcal{P}_F(tx+(1-t)y; \Omega_{i_0})$, one has
\begin{align*}
C_F(tx+(1-t)y;\Omega_{i_0})&=\rho_F(\omega-(tx+(1-t)y)\\
&=\rho_F(t(\omega-x)+(1-t)(\omega-y))\\
&\leq t\rho_F(\omega-x)+(1-t)\rho_F(\omega-y)\\
&\leq tC_F(x;\Omega_{i_0})+(1-t)C_F(y;\Omega_{i_0}).
\end{align*}
The equality (\ref{quali}) implies
\begin{equation*}
\rho_F(t(\omega-x)+(1-t)(\omega-y))=t\rho_F(\omega-x)+(1-t)\rho_F(\omega-y)=\rho_F(t(\omega-x))+\rho_F((1-t)(\omega-y)).
\end{equation*}
Since $x,y \notin \Omega_{i_0}$, one has $\omega-x,\;\omega-y\neq
0$, and hence, by Lemma \ref{ln}, there exists $\lambda>0$ such that
\begin{equation*}
t(\omega-x)=\lambda (1-t)(\omega-y).
\end{equation*}
This implies
\begin{equation*}
\omega-x=\gamma (\omega-y), \mbox{ where }\gamma:=\dfrac{\lambda(1-t)}{t}\neq
1.
\end{equation*}
Thus,
\begin{equation*}
\omega=\dfrac{1}{1-\gamma}x-\dfrac{\gamma}{1-\gamma}y\in L(x,y),
\end{equation*}
which is a contradiction. The proof is now complete. $\h$

\begin{Proposition} Let $F$ and $\Theta_j$ be strictly convex for every $j\in J$, where $J\neq\emptyset$, and let $S$ be convex. Suppose that for any $x, y\in S$, $x\neq y$, there exists $j\in J$ such that
\begin{equation*}
L(x,y)\cap \Theta_j=\emptyset.
\end{equation*}
Under the assumption that for every $x\in S$, the projection set
$\Pi_F(x;\Theta_i)\neq \emptyset$ for all $j\in J$, the function
\begin{equation*}
\mathcal{T}_2(x):=\sum_{j\in J} T_F(x;\Theta_j)
\end{equation*}
is strictly convex on $S$.
\end{Proposition}
{\bf Proof: }Suppose by contradiction that there exist $x\neq y$, $x,y\in S$, and $t\in (0,1)$ such that
\begin{equation*}
\mathcal{T}_2(tx+(1-t)y)=t\mathcal{T}_2(x)+(1-t)\mathcal{T}_2(y).
\end{equation*}
Using the convexity of each $T_F(x;\Theta_j)$ for $j\in J$, one has
\begin{equation}\label{eql}
T_F(tx+(1-t)y; \Theta_j)=tT_F(x;\Theta_j)+(1-t)T_F(y; \Theta_j) \mbox{ for all }j\in J.
\end{equation}
Suppose that $L(x,y)\cap \Theta_{j_0}=\emptyset$, where $j_0\in J$. Define
\begin{equation*}
u:=\Pi_F(x;\Theta_{j_0}) \mbox{ and } v:=\Pi_F(y;\Theta_{j_0}).
\end{equation*}
Then equation (\ref{eql}) implies
\begin{align*}
t\rho_F(u-x)+(1-t)\rho_F(v-y)&=tT_F(x;\Theta_{j_0})+(1-t)T_F(y; \Theta_{j_0})\\
&=T_F(tx+(1-t)y; \Theta_{j_0})\\
&\leq \rho_F(tu+(1-t)v-(tx+(1-t)y))\\
&\leq t\rho_F(u-x)+(1-t)\rho_F(v-y).
\end{align*}
It follows that $tu+(1-t)v=\Pi_F(tx+(1-t)y; \Theta_{j_0})$. This
implies $u=v$, since otherwise, $tu+(1-t)v\in \mbox{int
}\Theta_{j_0}$, which is a contradiction. Thus,
$u=v=\Pi_F(tx+(1-t)y; \Theta_{j_0})$. Equation (\ref{eql}), again,
implies
\begin{equation*}
\rho_F(u-(tx+(1-t)y)=\rho_F(t(u-x)+(1-t)(u-y))=\rho_F(t(u-x))+\rho_F((1-t)(u-y)).
\end{equation*}
Since $x,y \notin \Theta_{j_0}$, one has $u-x, u-y\neq 0$. Following
the proof of the previous proposition, one has
\begin{equation*}
u\in L(x,y),
\end{equation*}
which is a contradiction. The proof is now complete.
$\h$\vspace*{0.05in}

We are now ready to establish sufficient conditions for the optimization problem (\ref{hr}) to have at most one optimal solution.

\begin{Theorem} Let $F$ and $\Theta_j$ be strictly convex for every $j\in J$, and let $S$ be convex.
Suppose that for any $x, y\in S$, $x\neq y$, there exists $i\in I$
such that
\begin{equation*}
L(x,y)\cap \Omega_i=\emptyset
\end{equation*}
or there exists $j\in J$ such that
\begin{equation*}
L(x,y)\cap \Theta_j=\emptyset.
\end{equation*}
Under the assumption that for every $x\in S$, the projection sets
$\mathcal{P}_F(x;\Omega_i)\neq\emptyset$ and
\\$\Pi_F(x;\Theta_i)\neq \emptyset$ for all $i\in I$ and $j\in J$,
the function $\mathcal{H}(x)$ defined in {\rm(\ref{hr})} is strictly
convex, and the optimization problem {\rm (\ref{hr})} cannot have
more than one solution.
\end{Theorem}
\noindent {\bf Proof: }We have
\begin{equation*}
\mathcal{H}(x):=\mathcal{C}_2(x)+\mathcal{T}_2(x).
\end{equation*}
We only need to prove that $\mathcal{H}$ is strictly convex on $S$. However, this follows from the previous two propositions. The proof is now complete. $\h$ \vspace*{0.05in}

Finally, we state the related results for the Fermat-Torricelli problem generated by singletons.
\begin{Corollary}  Suppose $F$ is strictly convex and $S$ is convex. For any finite collection of points $\{a_i: i=1,\ldots,m\}$, $m\geq 1$, the optimization
\begin{equation*}
\mbox{\rm minimize }\sum_{i=1}^m\rho_F(a_i-x) \mbox{\rm \; subject to }x\in S
\end{equation*}
has at most one optimal solution, provided that $a_i$ for $i=1,\ldots,m$ are not collinear. This problem has a unique optimal solution if we assume additionally that $X$ is reflexive or $S$ is compact.
\end{Corollary}

\section{Concluding Remarks}

In this paper, we study generalized versions of the Sylvester
problem and the Fermat-Torricelli problem in Banach spaces. In the case where the space is $\Bbb R^n$ with the Euclidean norm, using generalized differentiation from convex
analysis, it is possible to construct explicitly solutions for
generalized Sylvester problems and generalized Fermat-Torricelli
problem for three arbitrary balls. Solutions for the
generalized Sylvester problems for three balls have a close
connection to the Apollonius' problem. These issues
are addressed in our paper \cite{nha}.

\end{document}